\documentclass[12pt]{amsart}

\usepackage{amssymb,url,mathrsfs,euscript}
\usepackage[spanish,british]{babel}

\theoremstyle{remark}     %roman body text
\newtheorem{rmk}{Remark}[section]\newtheorem*{rmk*}{Remark}

\newtheorem*{ex*}{Example}

\newtheorem*{acknowledgements}{Acknowledgements}

\swapnumbers              %places numbers before the heading
\theoremstyle{plain}      %default one (italics)

\newtheorem{proposition}[rmk]{Proposition}
\newtheorem{theorem}[rmk]{Theorem}

\theoremstyle{definition} %roman body

\newcommand{\be}{\begin{equation}}
\newcommand{\ee}{\end{equation}}
\newcommand{\bdm}{\begin{displaymath}}
\newcommand{\edm}{\end{displaymath}}

\numberwithin{equation}{section}

\newcommand{\hodge}{{*}}
\newcommand{\lie}[1]{\mathfrak{#1}}
\newcommand{\Lie}[1]{\textsl{#1}}
\newcommand{\G}{G\sb 2}
\DeclareMathOperator{\SO}{\Lie{SO}}
\DeclareMathOperator{\SL}{\Lie{SL}}
\DeclareMathOperator{\Spin}{\Lie{Spin}}
\DeclareMathOperator{\SU}{\Lie{SU}}

\DeclareMathOperator{\su}{\lie{su}}
\DeclareMathOperator{\so}{\lie{so}}

\newcommand{\g}{\lie{g}}
\newcommand{\n}{\lie{n}}
\newcommand{\s}{\lie{s}}

\renewcommand{\a}{\lie{a}}

\newcommand{\W}{\mathcal{W}}

\newcommand{\T}{\mathcal{T}}

\newcommand{\psp}{\psi^+}
\newcommand{\psm}{\psi^-}

\newcommand{\iso}{\cong}

\newcommand{\lto}{\longrightarrow}

\newcommand{\lan}{\langle}
\newcommand{\ran}{\rangle}
\newcommand{\hook}{\lrcorner\,}

\newcommand{\rcomp}[1]{\mathopen{\big[\mkern-5mu\big[} #1
\mathclose{\big]\mkern-5mu\big]}}
\newcommand{\rreal}[1]{\bigl[#1\bigr]}

\newcommand{\ol}{\overline}
\newcommand{\sym}{\mathcal{S}}

\newcommand{\oo}{\omega^2}
\newcommand{\1}{{\rm e}^1}\newcommand{\2}{{\rm e}^2}\newcommand{\3}{{\rm e}^3}
\newcommand{\4}{{\rm e}^4}\newcommand{\5}{{\rm e}^5}\newcommand{\6}{{\rm e}^6}
\newcommand{\7}{{\rm e}^7}
\renewcommand{\geq}{\geqslant}

\newcommand{\ad}{\textsl{ad}}
\newcommand{\tr}{\textsl{tr}}

\newcommand{\ds}{\displaystyle}\newcommand{\tx}{\textstyle}

\newcommand{\q}{\quad}\newcommand{\qq}{\qquad}

\newcommand{\bproof}{\begin{proof}}
\newcommand{\eproof}{\end{proof}}

\newcommand{\bsub}{\begin{subequations}}
     \newcommand{\esub}{\end{subequations}}
\newcommand{\ba}{\begin{array}}\newcommand{\ea}{\end{array}}

\newcommand{\f}{\varphi}\newcommand{\ff}{\hodge\varphi}
\newcommand{\R}{\mathbb{R}}\newcommand{\C}{\mathbb{C}}

\newcommand{\La}{\Lambda}

\newcommand{\w}{\wedge}
\newcommand{\set}{\subseteq}\newcommand{\sset}{\subset}
\newcommand{\na}{\nabla}

\newcommand{\vs}{\vphantom{\int_W^M}}
\newcommand{\VS}{\vphantom{\biggl\{ }}

\begin{document}
%%%%%%%%%%%%%%%%%%%%%%%%%%%%%%%%%
%%%THIS IS THE ARXIVE VERSION%%%%
%%%%%%%%%%%%%%%%%%%%%%%%%%%%%%%%%
\selectlanguage{british} 

\title[Special $\G$ metrics]{Special metrics in
$\boldsymbol{\G}$ geometry}
%-------------------------
\author[S.Chiossi]{Simon G.~Chiossi} \address[S.Chiossi]{Institut 
f\"ur Mathematik,
   Humboldt-Universit\"at zu Berlin, Unter den Linden 6,
   10099 Berlin, Germany}\email{sgc@math.hu-berlin.de}
%------------------------
\author[A.Fino]{Anna Fino} \address[A.Fino]{Dipartimento
      di Matematica, Universit\`a di Torino, via Carlo Alberto 10, 10123
      Torino, Italy} \email{fino@dm.unito.it}
%-----------------------
\begin{abstract}
We review a recent series of $\G$ manifolds
constructed via solvable Lie groups obtained in \cite{Chiossi-F:CG2}. 
They carry two related
distinguished metrics, one negative Einstein and
the other in the conformal class of a Ricci-flat metric.
\end{abstract}
%----------------------
\thanks{Supported by GNSAGA of INdAM and MIUR (Italy), and the SFB 647 
{\sl `Space -- Time -- Matter'} of the DFG}
%%%%%  \date{\today}
\maketitle
%---------------------

\section{Introduction}
\noindent 
A seven-dimensional Riemannian manifold is called a $G_2$-manifold
whenever the structure group of the tangent bundle is
contained in the  subgroup $G_2$ of the orthogonal group. 
Admitting such a reduction is
equivalent to the existence of a non-degenerate  three-form of positive type
$\varphi$. When this form is covariantly constant with respect to
the Levi--Civita connection then the holonomy of the manifold is
contained in $G_2$,
and the corresponding manifold is called parallel. 
Despite major advances in understanding exceptional geometry, 
producing metrics with holonomy equal
to $\G$ is still today not that easy. A quick browse through the literature of
high-energy physics gives evidence of this, although a rather good
assortment is available. Exhibiting 
\emph{complete} metrics, instead, has remained arduous
since the first examples
were constructed by Bryant and Salamon \cite{Bryant-S:exceptional}
almost twenty years ago: these are built from $\SO(4)\sset\G$
on vector bundles over $3$- and $4$-manifolds, 
and have been
relentlessly referred to ever since, hence becoming somehow
`classical'. We recall one construction in Section \ref{sec:oldmetrics}, 
by which the virtues of that landmark paper will be even more
apparent. 
\smallbreak

\noindent This note consists of two parts, the first of which pretentiously
tries to collect material on $\G$ structures placing
emphasis on Riemannian metrics. Needless to say, the
exposition is incomplete, hence we recommend to begin with 
 \cite{Bryant:exceptional}, \cite{Salamon:holonomy} and a recent 
survey \cite{Salamon:tour}.
The second half reviews the results obtained in a previous paper 
\cite{Chiossi-F:CG2} and 
describes in some detail two solutions
of the Hitchin flow given by metrics with 
holonomy equal to the full $G_2$. 
\smallbreak

Two incomplete holonomy metrics  with a $2$-step  nilpotent isometry
group $N$, whose orbits are hypersurfaces realised as torus bundles 
over tori, were presented in \cite{Gibbons-LPS:domain-walls}. It was
shown that such
Ricci-flat metrics are intimately related to complete
Einstein manifolds with a transitive solvable group of isometries. The metrics 
arise from Heisenberg limits of the isometry group 
of the two complete cohomogeneity-one metrics of
\cite{Bryant-S:exceptional}, once again leading back to these
examples.
They are moreover scale-invariant, that is to say they have additional symmetries 
generated by a homothetic Killing vector.
\smallbreak

In \cite{Chiossi-F:CG2} we showed that the previous Ricci-flat metrics are
 conformal to complete homogeneous metrics on a special kind of
solvable 
Lie group.
This is a rank-one solvable extension of the previous $N$, 
 now equipped with an $\SU(3)$-structure 
$(J, \omega, \psi^+)$, where $(J, \omega)$ is an almost
Hermitian structure and $\psi^+$ a $3$-form of unit norm.
The extension is constructed using a derivation $D$ on the metric Lie algebra
$(\n, \lan\,,\, \ran)$ of $N$. Such a
derivation is non-singular,  self-adjoint with 
respect to $\lan \,,\, \ran$ and satisfies $(DJ)^2 = (JD)^2$, a relation 
of primary importance for the construction.

Conformally parallel
$G_2$ structures on rank-one solvable extensions of $6$-dimens\-ion\-al
nilpotent Lie algebras $\n$ endowed with an $SU(3)$ structure and a
  derivation $D$ satisfying the previous requirements were then studied. 
We described thoroughly the corresponding
metrics with holonomy contained in $\G$ obtained 
after the conformal change. These fall into two categories: one where the 
holonomy is properly contained in the exceptional group, and the other consisting of
structures with holonomy truly equal to $\G$, depending upon whether 
the Lie algebra is irreducible. 
With this one is spared the sometimes daunting
task of detecting special spinor fields.
\smallbreak

As a consequence of the fact that $g$ belongs to the same conformal
class of the holonomy metric, the $\SU(3)$ structure on  the
 Lie group N is half-flat \cite{Chiossi-S:SUG}. 
A breakthrough technique
introduced by Hitchin \cite{Hitchin:forms} predicts then that there
is a (usually incomplete) metric $\tilde g$ with $\textsl{Hol}(\tilde
g)\set\G$ 
on the
product of $N$ with some real interval, but this is hard to determine
explicitly.
Nevertheless, due to the convenient set-up at hand, 
we computed the solution of the evolution
equations for each of the half-flat $\SU(3)$ structures. See
Section \ref{sec:examples} where full computations are carried out for the
structures arising from the nilpotent Lie groups isomorphic to
\smallbreak

$[{\rm e}_2,{\rm e}_1]={\rm e}_4,\q [{\rm e}_3,{\rm e}_1]={\rm e}_5,\q 
    [{\rm e}_3,{\rm e}_2]={\rm e}_6$ \q and \\[4pt]
\indent 
  $[{\rm e}_2,{\rm e}_1]={\rm e}_5,\q [{\rm e}_4,{\rm e}_1]=[{\rm e}_3,
{\rm e}_2]={\rm e}_6$.
\smallbreak

\noindent Comparing the 
solutions found in this way with the metric obtained with
the conformal change we proved that the metrics,  
though arising by two completely different methods, coincide, as one would expect.

As this were not enough to raise interest in this class of structures already, 
current work on these and other solvable examples shows
a peculiar spinorial behaviour \cite{Agricola-CF:strongweakG2}.
\medbreak
%-----------------------
\begin{acknowledgements}
This article was conceived 
to follow the {\sl II Workshop in Differential
  Geometry} in La Falda, Argentina. 
The authors are truly grateful to the organisers.

Untold thanks are due to R.Cleyton, S.Console and S.Garbiero for reading
the manuscript, and S.Salamon for comments.
\end{acknowledgements}

%%%%%%%%%%%%%%%%%%%%%%%%%%%%%%%%%%%%%%%%%%%%%%%%%%%%%
\section{$SU(3)$ and $G_2$ structures}\label{sec:sug}

\noindent 
Let $N^6$ be a real six-dimensional manifold.  An $SU(3)$-reduction
is given by an almost Hermitian triple $(J,\,h,\,\omega)$ where $h$ is 
a Riemannian metric, $J$ an $h$-orthogonal complex structure, $\omega$ the 
induced $(1,1)$-form, together with $\psi$ a $(3,0)$-form of unit
norm. Let 
\[\psi^+=\tfrac 12 (\psi+\ol{\psi}),\q \psi^-=\tfrac 1{2i}(\psi-\ol{\psi})\]
be the real and imaginary parts of $\psi$.
There is an orthonormal basis of 1-forms such  that
\be\label{eq:su3}
\begin{array}{c}
\omega={\rm e}^{14}-{\rm e}^{23}+{\rm e}^{56},\\[4pt]
\psi=({\rm e}^1\!+\!i{\rm e}^4)\wedge({\rm e}^2\!-\!i{\rm e}^3)\wedge
({\rm e}^5\!+\!i{\rm e}^6).
\end{array}
\ee
In the paper we will always indicate by $({\rm e}^{i})$ a -- possibly local -- 
basis of $1$-forms, 
so ${\rm e}^{ij}={\rm e}^{i}\w{\rm e}^{j}$. Dual vectors will be denoted by lower 
indexes, so ~${\rm e}^k({\rm e}_l)=\delta_{kl}$.
\smallbreak

Observe that the differential forms defining the reduction satisfy
 $\psi\wedge\omega=0$ and $3\psi\w\ol{\psi}=4i\,\omega^3$.
Since  $\psi^+$ is chosen to have stabiliser $\SL(3,\C)$ in the
general linear group, it determines the almost complex
structure $J$ and $\psi^-=J\psi^+$ \cite{Hitchin:forms}.

With $\su(3)^\perp=\so(6)/\su(3)$ denoting the orthogonal complement 
of $\su(3)$ in $\so(6)$, the known identifications
\bdm
\rreal{\Lambda^{1,1}_0}\iso\su(3),\qq 
\R\oplus\rcomp{\Lambda^{2,0}}\iso\su(3)^\perp %=\so(6)/\su(3)
\edm
allow one to split the space $T^*N\otimes\su(3)^\perp=\W$ in
$7$ irreducible $\SU(3)$-submodules
$$
\begin{array}{rcl}
   %    T^*N\otimes \su(3)^\perp 
  \W & = &
       \W_1^+\oplus\W_1^-\oplus\W_2^+\oplus\W_2^-\oplus\W_3\oplus\W_4\oplus\W_5.
     \end{array}
$$ 
The intrinsic torsion, a tensor $\tau \in \W$,
 accounts for the `non-integrable' structures, for the holonomy 
of the Riemannian metric $h$ is contained in $SU(3)$ if 
and only if $\omega, \psi^+$ and $\psi^-$ are all closed, in other
words when $\tau=0$.

The $\SU(3)$-representations $\W_i$ extend the original
  Gray--Hervella classes, and the reader might be familiar with some
names, like \smallbreak

- nearly K\"ahler structures, for which
$\tau\in\W_1^+\oplus\W_1^-$;  \smallbreak

- symplectic structures, where
$\tau\in\W_2^+\oplus\W_2^-\oplus\W_5$;  \smallbreak

- Hermitian structures, corresponding to  
$\tau\in\W_3\oplus\W_4$,  \smallbreak

\noindent but will probably not be acquainted with
the remarkable \smallbreak

- half-flat class $\W_1^-\oplus\W_2^-\oplus\W_3$.  \smallbreak

\noindent It is easily seen that
picking $\tau$ in this space is the same as demanding that 
$\psi^+$ and $\oo=\omega\w\omega$ be closed forms. The
name is designed to remind the fact that of the original dimension of
$\W$, only half survives.
\medbreak

We say that a $7$-dimensional manifold $M^7$ is {\sl built} from $N$ if the
cotangent space of $M$ splits at each point $m \in M$ as
\[T^*\sb m M^7 = T^*_nN^6\oplus
\R.\]
Hypersurfaces $N\sset M$, fibre bundles $M\lto
N$, or quotients $M/S^1$ are instances thereof. As $\SU(3)$ is a maximal
subgroup of $\G$, the special Hermitian geometry of $N$ 
induces a differential form on $M^7$ (pullbacks omitted)
\be\label{eq:sug}
\begin{array}{rcl}
\varphi &=&\omega\wedge {\rm e}^7+\psi^+,\\[4pt]
\end{array}
\ee
where ${\rm e}^7 $ is a $1$-form on $\R$. The
three-form $\f$ has isotropy 
$G_2\subset SO(7)$ and  determines a
compatible Riemannian metric $g$ and the
$4$-form $\ff=\psi^-\wedge {\rm e}^7+\frac 12 \,\omega^2
$, via the Hodge operator $*$. Implementing 
the basis of $T^*N$ with $\7$ preserves orthonormality and 
\be\label{eq:phi}
\hbox{$\varphi\!=\!
{\rm e}^{147}\!-\!{\rm e}^{237}\!+\!{\rm e}^{567}\!+\!{\rm
  e}^{125}\!+\!
{\rm e}^{136}\!+\!{\rm e}^{246}
\!-\!{\rm e}^{345}.$}
\ee
If (and only if) $M$ is parallel, 
$\varphi$ and $\ff$ become closed \cite{Fernandez-G:G2} and the induced 
metric $g$ has zero
Ricci curvature \cite{Bonan:G2-Spin7}.

The intrinsic torsion of a $G_2$ structure can be identified with
the covariant derivative  of the fundamental form
with respect to the Levi-Civita connection $\na$.
In \cite{Fernandez-G:G2} (see also \cite{Cabrera-MS:G2})
 a classification of $G_2$-manifolds in 16 
classes is given by studying the $G_2$-irreducible components of the 
torsion space $\T$.
Fern\'andez and Gray proved that $\T$ consists 
 of tensors having the same symmetries as $\nabla \varphi$  and
has four $G_2$-irreducible
components $\T_i$, $i = 0, \ldots, 3$.

On a $G_2$-manifold, the group's action on the tangent spaces
  $T_mM\iso \R^7=V$ induces an
  action on the exterior algebra $\Lambda^p (M)$. There are
  decompositions into modules
$$
\begin{array} {l}
\Lambda^2 V^* = \Lambda^2_{14} \oplus
  \Lambda^2_7 \iso \Lambda^5 V^*,\\[4pt]
\Lambda^3 V^* = \Lambda^3_{27}  \oplus 
\Lambda^3_7 \oplus \Lambda^3_1 \iso \Lambda^4 V^*,
\end{array}
$$
where $\Lambda^k_p $ denotes a certain irreducible $G_2$-module of
  dimension $p$.  The intrinsic torsion
of the $G_2$-structure  is encoded in the exterior 
derivatives $d\f, d\ff$ as follows
$$
\begin{array}{l}
d \varphi = \tau_0 \ff + 3 \tau_1 \wedge \varphi + 
\hodge\tau_3,\\
d \ff = 4 \tau_1 \wedge \ff + \tau_2 \wedge \varphi,
\end{array}
$$
for unique differential forms 

$\tau_0 \in \R \cong \T_0$,\\[2pt]
\indent $\tau_1 \in \La^1 \cong \T_1$, \\[2pt]
\indent $\tau_2 \in 
\Lambda^2_{14} \cong\g_2\ \text{(the exceptional Lie algebra)}\ 
\cong \T_2$, \\[2pt]
\indent $\tau_3 \in \Lambda^3_{27} \cong \sym^2_0 V^*\ 
\text{(the space of traceless sym\-met\-ric $2$-tensors)}\ \cong \T_3$, \\[2pt]
\noindent see for instance \cite{Cabrera:G2,Bryant:G2}.
\smallbreak

Friedrich and Ivanov proved that  $\tau_2 =0$
if and only if there exists an affine connection $\tilde \nabla$ with
totally skew-symmetric torsion $T$ such that $\tilde \nabla \varphi =
0$ \cite{Friedrich-I:skew}.
Then $M\sp7$ is a `$G\sb2$-manifold with torsion' ($\G$T) and the resulting 
torsion $3$-tensor is 
$$
T = \tfrac 76 \tau_0 \varphi -* d\varphi +
*(4\tau_1\wedge\varphi).$$
\smallbreak

\noindent An interesting subset of $G_2$T-manifolds consists of
those of type $\T_1$,  for which
$$
\begin{array} {c}
d \varphi = 3 \tau_1 \w\f,\\
d (\ff) = 4 \tau_1 \w\ff
\end{array}
$$
so
\bdm
T = \hodge(\tau_1 \wedge \varphi).
\edm
These manifolds are also called locally conformally parallel, since 
the change
$e^{2f} g$  (with $df = - \frac {1} {3} \tau_1$) gives locally a 
parallel structure. % a metric with holonomy contained in $G_2$.
\medbreak

The reader interested in {\sl compact} $\G$ manifolds of class $\T_1$
should look at \cite{Ivanov-PP:lcp-exceptional}.

%%%%%%%%%%%%%%%%%%%%%%%%%%%%%%%%%%%%%%%%%%%%%%%%%%%%%%%%%%%%%
\section{Some examples of $\G$ metrics}\label{sec:oldmetrics}

\noindent 
{\bf 1.} We recall one essential idea of \cite{Salamon:exceptional}.
Let $K$ be an oriented Riemannian $4$-manifold with local
orthonormal basis ${\rm f}^4,{\rm f}^5,{\rm f}^6,{\rm f}^7$ of the
cotangent bundle. Define the unit forms
\be\label{eq:sd}
{\rm f}^1={\rm f}^4{\rm f}^5-{\rm f}^6{\rm f}^7,\q 
{\rm f}^2={\rm f}^4{\rm f}^6-{\rm f}^7{\rm f}^5,\q 
{\rm f}^3={\rm f}^4{\rm f}^7-{\rm f}^5{\rm f}^6
\ee
to span $\La^2_-T^*K$. The total space $Y$ of the
latter decomposes as $H\oplus V$. The vertical space $V$ is generated by
three $1$-forms $({\rm e}^j)$  on $Y$ depending on the fibre coordinates,
whilst $H$ has basis (the pullbacks of)
\eqref{eq:sd}. Given now two positive functions $\alpha,\beta$ on
$Y$,
\bdm
\phi=6\alpha^3\,{\rm e}^1{\rm e}^2{\rm e}^3+
\alpha\beta^2\, d({\rm f}^1{\rm e}^1+{\rm f}^2{\rm e}^2+{\rm f}^3{\rm e}^3)
%%% d(\tx\sum {\rm f}^j{\rm e}^j)
\edm
is a $\G$ structure determining a Riemannian metric of the form
$\alpha^2g_V+\beta^2g_H$ in terms of the above splitting. Now if $K$
is self-dual and positive Einstein, 
choosing $\beta=(tr+1)^{1/4},\ \alpha=\beta^{-1}$, with 
$r>0$ a radial coordinate and some positive $t$, renders $\f$
closed and coclosed, hence parallel, and the metric 
\be\label{eq:BS}
(tr+1)^{-1/2}g_V+(tr+1)^{1/2}g_H
\ee
is complete, Ricci-flat and has holonomy equal to $\G$. 
When the parameter tends to zero, the metric becomes
conical on the product of $\R^+$ and the twistor space.
Since $K$ is 
$S^4$ or $\mathbb{CP}^2$, the groups
$\SO(5)$, $\SU(3)$ act isometrically with generic orbits of
codimension one. 

The metric resembles
the Eguchi-Hanson instanton \cite{Eguchi-H:self-dual}, which is
Einstein on $T^*\mathbb{CP}^1$ and makes the standard
holomorphic symplectic form covariantly constant.
\smallbreak

{\bf 2.} A similar example, constructed in the
flavour of Section \ref{sec:sug}, is the following. 
Let $Iw$ be the compact quotient
of the complex $3$-dimensional Heisenberg group by Gaussian
integers, called Iwasawa manifold. The product $Iw\times \R$ admits an orthonormal
basis $({\rm e}^j)$ with ${\rm e}^j={\rm f}^j$ of \eqref{eq:sd} for
$j\geq 4$, such that
\bdm
d{\rm e}^j=\left\{
\ba{ll} 
{\rm f}^j & j=1,2,3\\
0 & j=4,5,6,7
\ea\right.
\edm
 and the three-form  
 \bdm
 \phi={\rm e}^{127}+{\rm e}^{347}+{\rm e}^{657}-
 \rm{e}^{135}+\rm{e}^{126}+\rm{e}^{643}+\rm{e}^{254}
 \edm
 is a $\G$ structure. Indicating by $\hook$ the interior product 
of a vector with a differential form, the invariant tensors on $Iw$
 \bdm
 \psi^+=d(\rm{e}^{56}),\q \omega={\rm e}_7\hook\phi
 \edm
 are such that $d\hodge\omega=0$ and $d\omega=\psi^+$. 
It is no coincidence that this almost complex structure recalls the one  
investigated in 
 \cite{Abbena-GS:aH-nil} as a distinguished element in a `twistor space' for $Iw$.
The $\SU(3)$ reduction $(\omega,\psi^+)$ is merely a modification of
\eqref{eq:su3} obtained by rotations in the bundle $Iw\lto T^4$
fibred by $2$-tori, reminding of the Penrose
fibration $\mathbb{CP}^3\stackrel{S^2}{\lto}S^4$. 
A more systematic approach including this example was developed
in \cite{Apostolov-S:K-G2}.
\smallbreak

{\bf 3.} A central chapter of the theory of 
exceptional geometry is related to  Killing
spinors \cite{Baum-FGK:spinors}.
This notion allowed B\"ar to prove that  \cite{Baer:real}
if $(X^6, g)$ is nearly K\"ahler, then the metric cone $(X^6 \times \R^+, 
t^2\,g+dt^2)$ has holonomy $\G$.
\smallbreak

{\bf 4.} 
Physical evidence has now shifted most of the concern 
towards metrics with orbifold singularities, see
\cite{Acharya-G:singularities-report,Acharya-W:chiral,Atiyah-W:G2}. 
%%% The simplest one is the following. 
One with an isolated conical singularity (the most subtle of the three
known in the simply-connected case) is the following. 
The space $X=S^3\times S^3$ admits Einstein metrics, the easiest being the
product of the two round metrics on the factors, that has symmetry
$\SU(2)\times\SU(2)\times\SU(2)\times\SU(2)$. It has another 
Einstein -- and here more relevant one invariant under
$\SU(2)^3\times\Sigma_3$, where the latter is the symmetric group 
on $3$ elements generating `triality'. 
Describing $X$ as the $3$-symmetric space
$\SU(2)^3/\SU(2)$ under a diagonal action, the metric is 
\bdm
g_X=-\tr \bigl( (a^{-1}da)^2+(b^{-1}db)^2+(c^{-1}dc)^2\bigr),
\edm
where $a=g_2g_3^{-1},\ b=g_3g_1^{-1},\ c=g_1g_2^{-1}$ and
$(g_1,g_2,g_3)\in \SU(2)^3$. The cone of $g_X$ deforms to a smooth
complete holonomy metric on some $Y$ \cite{Bryant-S:exceptional}, 
itself homeomorphic to $\R^4\times S^3$, because in the limit one
of the spheres $S^3$ collapses. Thus $Y$ has an asymptotically conical $\G$-metric.
\smallbreak

{\bf 5.} The striking results achieved with the discovery of compact
manifolds with holonomy $G_2$ by Joyce 
\cite{Joyce:G2-1&2} first, and Kovalev 
\cite{Kovalev:G2}  by different methods, answered the 
$\G$-analogue of the Calabi conjecture on special Hermitian holonomy. This is the
origin of the expression {\sl Joyce manifolds}.
These constructions do not yield explicit metrics, though it must be
said that often they need not be, at least for the purposes of
string theorists.

%%%%%%%%%%%%%%%%%%%%%%%%%%%%%%%%%%%%%%%%%%%%%%%%%%%%%%%
\section{Solvable extensions of nilpotent Lie algebras}
\label{sec:extensions}

\noindent 
  Let $N$ be now a 6-dimensional one-connected real nilpotent Lie
group. 
It is nilpotent if and only if
there exists a basis $({\rm e}^1, \ldots, {\rm e}^6)$ of
left-invariant 1-forms on $N$ such that
$$
d {\rm e}^i \in \Lambda^2 \left<{\rm e}^1, \ldots, {\rm e}^{i -
  1}\right>,\ i = 1, \ldots, 6.
$$
In terms of the lower central series $\n^0 = \n,\ \n^i
= [\n^{i - 1}, \n]$, this is the same as requiring $\n^s = 0$ for some 
$s\in{\mathbb N}$. 
If $N$ has  rational structure constants, then by
\cite{Malcev:rational} it admits a compact quotient $\Gamma\backslash N$ 
by a uniform 
discrete subgroup. Such a homogeneous space is called  nilmanifold.
\smallbreak

Solvable extensions of nilpotent Lie groups are particular examples of 
homogeneous Einstein spaces of negative scalar curvature. On the other
hand all known non-compact, non-flat, homogeneous Einstein
spaces have the form $(S, g)$, where $S$ is a solvable Lie group and
$g$ is a left-invariant metric, which we will denote by the name 
\emph{solvmanifold}.
Because left-invariant Einstein metrics 
on unimodular solvable Lie groups are flat
\cite{Dotti:unimodular-solvable}, 
the solvable Lie groups  we consider
will be not unimodular, hence never admit a compact quotient 
\cite{Milnor:curvature}.

The Einstein solvmanifolds available as of today  are modelled on 
completely solvable Lie groups -- the eigenvalues of
$ad_U$ are real, for any vector $U$ -- and their underlying metric Lie 
algebras $\big(\s,\,\lan\,, \ran\big)$ are {\sl standard} and {\sl of Iwasawa type}.
Given a metric nilpotent Lie algebra $\big({\n}, [\, 
,  ]_\n, \lan \, , \, \ran'\big)$ with inner
product $\lan\, ,\, \ran'$, a metric solvable Lie algebra
$\big({\s} = {\n} \oplus {\a}, [\,  ,\,  ], \lan \, , 
\, \ran\big)$ is called a metric solvable extension of  $\big({\n}, 
\lan \, ,\, \ran'\big)$ if $[\,  ,\,  ]$ restricted
to $\n$ coincides with $[\,  , \, ]_\n $ and
$\langle\, ,  \rangle \vert_{\n \times \n} = \langle\, ,  \rangle'$.
One says that $\s$ is standard if $\a =[\s,\s]^{\perp}$ is
Abelian.  The dimension of $\a$ 
is called the algebraic rank of $\s$.

\noindent 
If the rank is one, say $\a=\ <\!A\!>$, the
extension is of Iwasawa type if

(i)\q $\ad_A\ne0$ is self-adjoint
with respect to $\langle\, ,\,\rangle$, and

(ii)\q $(\ad_A)\vert_\n$ is positive-definite.

\noindent By
\cite[4.18]{Heber:noncompact-Einstein}  the study of standard 
Einstein  metric solvable Lie algebras  reduces to
rank-one metric solvable extensions
$$
\big(\s = \n \oplus \R H,\ \langle \cdot, \cdot \rangle\big)
$$
for some $\n$ and $H$ with $\langle H, \n
\rangle = 0$, $\vert \vert  H  \vert
\vert = 1$. The extended Lie bracket follows the rule
$$
\left\{\begin{array}{l} [H, X] = D (X),\\[4pt] [X, Y] = [X,
Y]_\n\end{array}\right.
$$
 and $D \in
{\mbox {Der}} (\n)$ is a  derivation of the Lie algebra.
If the metric $\langle\, ,\,\rangle$ on the extension $\s$ is Einstein,
then the derivation $D$ is necessarily 
 self-adjoint, and in fact unique \cite{Heber:noncompact-Einstein}.

For an Einstein solvmanifold, Heber calls eigenvalue
 type the sequence 
$(\lambda_1 < \ldots < \lambda_r;\ m_1, \ldots, m_r)$, where $(\lambda_i)$ are 
the eigenvalues of 
$D$ and $(m_i)$ the corresponding multiplicities. He proved that in any dimension
 only finitely many eigenvalue types
 occur. In addition, six is a critical dimension, since by 
\cite{Lauret:variational-method, Will:Einstein-7-solv} all 
nilpotent Lie groups up to dimension $6$ admit an extension
of rank one carrying an Einstein metric.

%%%%%%%%%%%%%%%%%%%%%%%%%%%%%%%%
\section{The homogeneous models}

\noindent 
From now $Y$ will indicate a manifold equipped with the conformally 
parallel $G_2$ structure \eqref{eq:phi}, so that
the holonomy group of the metric $e^{2f} g$ is contained in $G_2$. 
The function $f$ is prescribed by $d f = -m {\rm e}^7, m\in\R^-$.

In order to use the underlying almost Hermitian geometry, we suppose that
$Y$ arises from a rank-one solvable extension $\s$ of a metric
nilpotent Lie algebra $\n$, whose Lie group $N$ is
endowed with an invariant $\SU(3)$
structure $(J,\omega, \psi^+)$ and a non-singular self-adjoint
derivation $D$, as in the Einstein case. 
In this way the algebraic structure of $\s$ 
blends in with the Riemannian aspects of $Y$.
We  require in fact that $(DJ)^2 =
(JD)^2$, a condition that translates into nice features of string
models \cite{Gibbons-LPS:domain-walls}.
Concretely, the solvable structure is defined by the nilpotent
$\n$ and by taking
\bdm
%\label{eq:D}
      D=\ad_{{\rm e}_7}.
\edm
A classification result establishes that $N$ cannot be arbitrary, 
even $p$-step nilpotent with $p=1$ or $2$. 
Under the above assumptions in fact, we proved that
\begin{theorem}\cite{Chiossi-F:CG2}
\label{thm:classification}
Let $\s=\n\oplus \R {\rm e}_7$ be a rank-one solvable
     extension determined by $\ad_{{\rm e}_7}$. Then the $G_2$ structure
     $\f=\omega\wedge\7+\psi^+$ defined on $\s$ is conformally
     parallel 
     if and only if $\n$ is isomorphic to one of:
     \begin{center}
       $(0,0,{\rm e}^{15},{\rm e}^{25},0,{\rm e}^{12})$, \q 
       $(0,{\rm e}^{54},{\rm e}^{64},0,0,0)$,
       \q $(0,{\rm e}^{54},{\rm e}^{15}+{\rm e}^{64},0,0,0,0)$, \\[4pt]
       $(0,0,e^{15} + e^{64},0,0,0)$,\q
       $(0,{\rm e}^{61}+{\rm e}^{54},{\rm e}^{15}+{\rm e}^{64},0,0,0)$, 
       \q $(0,0,{\rm e}^{15},0,0,0)$,\\[4pt]
        $(0,0,0,0,0,0)$.
     \end{center}
   \end{theorem}
\noindent Explicitly, the Lie algebras
found are listed in the Table that follows.
The terms corresponding to the nilpotent part have been
highlighted to make it easier to recognize the underlying
$\n$ of Theorem \ref{thm:classification}. \\
About the notation: the 
`differential' expression
$(0,0,{\rm e}^{15},{\rm e}^{25},0,{\rm e}^{12},0)$ is a quick way of saying
$[{\rm e}_5,{\rm e}_1]={\rm e}_3,\ [{\rm e}_5,{\rm e}_2]={\rm e}_4,
\ [{\rm e}_2,{\rm e}_1]={\rm e}_6$ for the basis of $\s$. In fact a general  
Lie algebra $\g$ of dimension $n$ is either prescribed by a bracket $[\, , ]$, 
or by a map $d:\g^*\lto \Lambda^2\g^*$ which 
extends to give a complex
\bdm
0\lto\g^*\lto\Lambda^2\g^*\lto \Lambda^3\g^*\lto\ldots\lto\Lambda^n\g^*\lto 0.
\edm

{\footnotesize
  \begin{table}
$$  \ba{|c|} 
      \hline
       \hbox{\normalsize Solvable Lie algebras}\VS \\
       \hline\hline
        (-m{\rm e}^{17},-m{\rm e}^{27},-m{\rm e}^{37},
	  -m{\rm e}^{47},-m{\rm e}^{57},-m{\rm e}^{67},0)\VS \\
       \hline
        (-\tfrac 23 m{\rm e}^{17},-m{\rm e}^{27},-\tfrac 43 
       m{\rm e}^{37}+\boldsymbol{\tfrac 23
        m{\rm e}^{15}}, -m{\rm e}^{47}, -\tfrac 23 m{\rm e}^{57},
-m{\rm e}^{67},0)\VS \\
\hline
        \bigl(-\tfrac 34 m{\rm e}^{17},-m{\rm e}^{27},
       -\tfrac 32 m{\rm e}^{37}\boldsymbol{+\tfrac 12 m({\rm e}^{15}-{\rm e}^{46})},
       -\tfrac 34 m{\rm e}^{47},
       -\tfrac 34 m{\rm e}^{57},-\tfrac 34 m{\rm e}^{67},0\bigr)\VS \\
       \hline
        \bigl(-\tfrac 45 m{\rm e}^{17},-\tfrac 65 
m{\rm e}^{27}\boldsymbol{-\tfrac 25 m{\rm e}^{45}},
        -\tfrac 75 m{\rm e}^{37}\boldsymbol{+\tfrac 25 
m({\rm e}^{15}-{\rm e}^{46})},-\tfrac 35 m{\rm e}^{47},
        -\tfrac 35 m{\rm e}^{57},-\tfrac 45 m{\rm e}^{67},0\bigr) \VS  \\
       \hline
        (-m {\rm e}^{17},-\tfrac 54 m{\rm
       e}^{27}\boldsymbol{-\tfrac 12 m{\rm e}^{45}},
        -\tfrac 54 m{\rm e}^{37}\boldsymbol{-\tfrac 12 m{\rm e}^{46}}, -\tfrac 12 
m{\rm e}^{47},-\tfrac
        34 m{\rm e}^{57}, -\tfrac 34 m{\rm e}^{67},0)\VS \\
       \hline
          \bigl(-\tfrac 23 m{\rm e}^{17},-\tfrac 43 
m{\rm e}^{27}\boldsymbol{-\tfrac 13 m({\rm e}^{16}+{\rm e}^{45})},
        -\tfrac 43 m{\rm e}^{37}\boldsymbol{+\tfrac 13 
m({\rm e}^{15}-{\rm e}^{46})},-\tfrac 23 m{\rm e}^{47},
        -\tfrac 23 m{\rm e}^{57},-\tfrac 23 m{\rm e}^{67},0\bigr)\VS \\
       \hline
       ( -\tfrac 35 m{\rm e}^{17},-\tfrac 35 m{\rm e}^{27},-\tfrac 65 
m{\rm e}^{37}\boldsymbol{+\tfrac
        25 m{\rm e}^{15}},
        -\tfrac 65 m{\rm e}^{47}\boldsymbol{+\tfrac 25 m{\rm e}^{25}},
-\tfrac 35 m{\rm e}^{57},
        -\tfrac 65 m{\rm e}^{67}\boldsymbol{+\tfrac 25 m{\rm e}^{12}},0)\VS \\
       \hline
       \ea $$
     \end{table}
}  %%%%% this ends the footnotesize of the table %%%%%%

The construction of Section \ref{sec:extensions} is particularly interesting since
the results of \cite{Heber:noncompact-Einstein,Will:Einstein-7-solv} 
ensure that
$Y$ will admit a homogeneous Einstein metric with negative scalar
curvature, and moreover a unique one if one chooses
the eigenvalues of $\ad_{{\rm e}_7}$.
\subsection{Example of an Einstein metric.}
%\begin{ex}
To end this section, we provide one of the Einstein metrics.
Consider the  $3$-step solvable Lie algebra  with structure equations
\be\label{eq:3step}
\begin{array} {l}
d {\rm e}^1 =  2 b\, {\rm e}^{17} + \sqrt{6} b\, {\rm e}^{26}\\
d {\rm e}^i = b\, {\rm e}^{i7}, \quad  i = 2, 4, 6\\
d {\rm e}^3 =  2 b\, {\rm e}^{37} - \sqrt{6} b\, {\rm e}^{46}\\
d {\rm e}^5 = 2 b\, {\rm e}^{57}-\sqrt{6} b\, {\rm e}^{24}\\
d {\rm e}^7 = 0,
\end{array}
\ee
where $b$ is real and not zero.
This is the last one in the Table in disguise, endowed 
with the 3-form \eqref{eq:phi}.
The $\G$ structure satisfies the conditions 
$$
\begin{array} {l}
d \varphi =  5 b\, ( {\rm e}^{1257} + {\rm e}^{1367} - {\rm e}^{3457}) - 
3 b(\sqrt{6} -1) {\rm e}^{2467},\\[3pt]
d (\ff) =  \sqrt{6}\, b(\sqrt{6} +1) {\rm e}^{23567} + \sqrt{6}\, b(\sqrt{6} -1)
({\rm e}^{12347} - {\rm e}^{14567}),
\end{array}
$$
% with
% $$
% \ff=  - {\rm e}^{1357} + {\rm e}^{1267}  - {\rm e}^{2457} - {\rm e}^{3467} -
% {\rm e}^{1234} - {\rm e}^{2356}  + {\rm e}^{1456}.
% $$
so it belongs to the class
${\mathcal T}_1 \oplus {\mathcal T}_3$ and the associated metric
$\sum_{i = 1}^7 ({\rm e}^i)^2$  is
Einstein with Ricci tensor ${\mbox {\sl Ric}} (g) = - 15 b^2 g$. We shall 
return to this example at the very end of this survey.
%\end{ex}

%%%%%%%%%%%%%%%%%%%%%%%%%%%
\section{Ricci-flat metrics}

\noindent Another aspect of the picture is that the almost Hermitian manifold $N$
is half-flat~\cite{Chiossi-S:SUG}. 
If one considers a $6$-manifold $N$ equipped with a reduction 
 $(\omega, \psi^+)$ that depends on a real parameter $t\in I$, let's say a
 `time-depending' $\SU(3)$-structure, then
 $N^6 \times I$ is a warped $\G$ manifold with fundamental form
 \bdm
\f=\omega\w dt+\psp.
\edm
 % \begin{array}{rcl}
 % d \varphi &=& (\hat d \omega - \frac {\partial \psi^+} {\partial t} )
 % \wedge dt + \hat d \psi^+.
 % ,\\[4pt]
 % d \ff &=& (\hat d \psi^- + \omega \wedge \frac {\partial \omega}
 % {\partial t}) \wedge dt + \omega \wedge \hat d \omega,
 % \end{array}$$
 % $where $\hat d$ denotes exterior differentiation on $N^6$.

 \noindent If $(N^6 \times I, \varphi)$ is parallel, the forms $(\omega^2, \psi^+)$
 evolve according to differential
equations  
 \begin{equation}\label{eq:flow}
   d\omega=\tfrac {\partial \psi^+}{\partial t},\qq 
   d(J\psp)=-\tfrac{\partial}{\partial t}(\tfrac 12 \oo),
 \end{equation}
 coming from the Hamiltonian flow of a functional.
 The opposite is also true. If $N$ is compact and \eqref{eq:flow} are 
 satisfied by closed forms $\psi^+$  and $\omega^2$ of suitable algebraic type, 
 there exists a metric with holonomy contained in $G_2$ on
 the product of $N^6$ with some interval $I$ \cite{Hitchin:forms}.
 
 The system \eqref{eq:flow} is tough to solve in general. To apply Hitchin's theorem,
 instead of considering a nilpotent Lie group $N$ we work with the
 associated nilmanifold and use the left-invariance of the forms.
 The Ricci-flat metrics thus found were described in
 \cite{Chiossi-F:CG2} and seen to coincide with homogeneous metrics possessing a
 homothetic Killing field. This is attained by comparing 
 the expressions, and bearing in mind that the  simply-connected solvable Lie group $S$
 (corresponding to $\s$) is diffeomorphic to $\R^7$, hence admits
 global coordinates of type $(x_1,\ldots,x_6,t)$. At the same time the
 metric can be seen as living on the product $\R\times
 \Gamma\backslash N$, where the nilmanifold $\Gamma\backslash N$ is
 $2$-step nilpotent, or $T^6$. This explains the bundle
 structure appearing, since \cite{Palais-S:2-step-nils} torus fibrations 
 over tori are, essentially, nilmanifolds of step-length two. The isometry
 between the two metrics is given by an appropriate choice of frame $({\rm e}^i)
 $, whence 
 \be\label{eq:c-metric}
 g=e^{-2mt}\tx\sum_{i=1}^7 ({\rm e}^i)^2
 \ee
 has holonomy a  subgroup of $\G$.

%%%%%%%%%%%%%%%%%%%%%%%%%%%%%%%%%%%%%%%%
\section{Two examples with full holonomy
}\label{sec:examples}

\noindent 
We carry out some calculations showing how the $\G$-holonomy
metrics are related to the Ricci-flat ones.

\subsection{First example.} It is clear that the Lie algebra
%\be\label{12,13+24}
{\small
\begin{displaymath}
   \bigl(-\tfrac 45 m{\rm e}^{17},-\tfrac 65 m{\rm e}^{27}-\tfrac 25 m{\rm e}^{45},
   -\tfrac 75 m{\rm e}^{37}+\tfrac 25 m({\rm e}^{15}-{\rm e}^{46}),
-\tfrac 35 m{\rm e}^{47},
   -\tfrac 35 m{\rm e}^{57},-\tfrac 45 m{\rm e}^{67},0\bigr)
\end{displaymath}
}
%\ee
extends the nilpotent Lie algebra with non-zero brackets
      \begin{displaymath}
        {\rm e}_2=[{\rm e}_5,{\rm e}_4], \q {\rm e}_3=[{\rm e}_6,{\rm
	    e}_4]=
[{\rm e}_1,{\rm e}_5].
      \end{displaymath}
Consider the Riemannian metric 
      \begin{multline}
        \label{eq:G2metric}
        g =   e^{-2m t} dt^2 + e^{-\tfrac 25 m t}(dx_1^2 + dx_6^2) +
        e^{-\tfrac 45 m t} (dx_4^2 + dx_5^2) + \\[3pt]
        \tfrac 9{25}m^2e ^{\tfrac 45 mt}(dx_3-\tfrac 23 x_1dx_5 + 
	\tfrac 23 x_4
        dx_6)^2 + \tfrac 9{25}m^2e^{\tfrac 25 mt}(dx_2+\tfrac23 x_4dx_5)^2,
      \end{multline}
      whose holonomy group is precisely $G_2$. 
\begin{proposition}
This $g$ is 
locally isometric to
\begin{multline*}
     \label{eq:isometricnewmetric}
     \phantom{MMM} ds^2= V^3dw^2 + V(du_1^2 + du_4^2) +
                         V^2(du_2^2 + du_3^2) + \\[3pt]
                         V^{-2}\bigl(dy_2+k(u_2du_4-u_1du_3)\bigr)^2+
                         V^2(dy_1+k u_2du_3)^2,\phantom{MMM}
\end{multline*}
on the product of $\R\times \mathcal{B}$ where $\mathcal{B}$ 
is the total space of a torus bundle
     \begin{displaymath}
      \mathcal{B}\lto T^4.
      \end{displaymath}
\end{proposition}
\begin{proof}
By demanding that
\begin{equation}\label{eq:coords}
        \begin{array} {c}
          {\rm e}^1 = e^{\tfrac 45 m t} dx_1,\q  
	  {\rm e}^2 = -\tfrac 35 m\,e^{\tfrac 65 m t}( dx_2 + \tfrac
	  23 x_4 dx_5),\\[3pt] 
	  {\rm e}^3 = -\tfrac 35 m\,e^{\tfrac 75 m t}( dx_3 - \tfrac 23  x_1 dx_5
          + \tfrac 23  x_4 dx_6),\q {\rm e}^4 = e^{\tfrac 35 m t}
	  dx_4, \\[3pt]
	  {\rm e}^5 = e^{\tfrac 35 m 
	    t} dx_5,\q {\rm e}^6 = e^{\tfrac 45 m t}dx_6,\q {\rm e}^7 = dt
        \end{array}
\end{equation}
be an orthonormal frame for $\n\oplus\R$, we have that 
$g$ as in \eqref{eq:c-metric} recovers
\eqref{eq:G2metric}. The local equivalence is established once we
indicate by $(u_i)$ the coordinates on $T^4$, by $y_1,y_2$ those 
on the fibres $T^2$, and have $w$ describe the seventh direction, with $V=m\,w$.
\end{proof}
\bigbreak

\noindent
\begin{theorem}
The family of Ricci-flat metrics arising from the $\SU(3)$ structures
solutions of \eqref{eq:flow} essentially coincides with \eqref{eq:G2metric}. 
\end{theorem}
\begin{proof}
We deform the starting $\SU(3)$ reduction
\begin{equation}
  \tfrac 12 \omega_0\w\omega_0=-{\rm e}^{2356}-{\rm e}^{1423}+{\rm e}^{1456},\q
  \psp_0={\rm e}^{125}-{\rm e}^{345}+{\rm e}^{136}+{\rm e}^{246}
\end{equation}
determined by \eqref{eq:su3} by forms
on $\n$ representing zero cohomology, so that
\bdm
\bigl(\,[\omega_0^2],[\psp_0]\,\bigr)=\bigl(\,[\omega(t)^2],[\psp(t)]\,\bigr) 
\q \text{in}\ H^4(\Gamma\backslash N,\R)\times H^3(\Gamma\backslash
N,\R)\q \text{for all $t$'s}.
\edm 
The four-form $\omega_0\w\omega_0$ flows under \eqref{eq:flow}
according to
\begin{displaymath}
\tfrac 12 \oo(t)=
-P(t)({\rm e}^{1423}+{\rm e}^{2356})+\bigl(P(t)+D(t)\bigr){\rm e}^{1456}
  \end{displaymath}
for smooth maps $P,D$ with $P(0)=1, D(0)=0$.

In the same way the three-form turns out to be
\begin{displaymath}
      \psp(t) = \bigl(M(t)+1\bigr)({\rm e}^{125}-{\rm e}^{345}+
	  {\rm e}^{246})+{\rm e}^{136},\qq M(0)=1.
\end{displaymath}
This almost gives the K\"ahler form as
\begin{displaymath}
   \omega(t)=\sqrt{P+D}\,({\rm e}^{14}+{\rm e}^{56})+\frac {5M'}{2m}\,{\rm e}^{23},
\end{displaymath}
with the dash denoting derivatives with respect to $t$. 
Notice how the expression respects the bundle structure of $\mathcal{B}$.
At this point it is anybody's guess to solve \eqref{eq:flow}, because
one does not know $\psm(t)$, which normally makes the system extremely
hard to tackle. To this end we define an orthonormal basis $(\lambda^{a_i}{\rm e}^{i})$, 
in which the choice of exponents $a_i$ ought to reflect the form of $\psp(t)$ above.
Picking 
\begin{displaymath}
a_i=(-1,1,2,-2,-2,-1),
\end{displaymath}
for example, one has
\begin{displaymath}
  \begin{array} {l}
    \psp=\lambda^{-2}({\rm e}^{125}-{\rm e}^{345}+
	  {\rm e}^{246})+ {\rm e}^{136},\\[3pt]
    \psm=\lambda^{-2}{\rm e}^{126}-\lambda^{-4}\,{\rm e}^{245}-\lambda^{-3}\,{\rm
      e}^{135}-\lambda^{-1}\,{\rm e}^{346}
  \end{array}
\end{displaymath}
giving
\bdm
% \label{eq:dpsm}
P=1,\q D'=\tfrac 65 m\sqrt{M+1}
\edm
with 
\bdm 
\lambda(t)^{-2}=M+1.
\edm
Since $\psp\w\psm$
and $\omega^3$ are both volume forms on $N$, the 
Cauchy system is solved by
\begin{displaymath}
M(t)=(1-mt)^{2/5}-1.
\end{displaymath}
It is then a simple matter to write the induced metric $g(t)$
      \begin{multline*}
       (1-mt)^{2/5}\,\bigl((\1)^2+(\6)^2\bigr)+
(1-mt)^{4/5}\,\bigl((\4)^2+(\5)^2\bigr)\\[3pt]
+(1-mt)^{-2/5}\,(\2)^2+(1-mt)^{-4/5}\,(\3)^2+dt^2
      \end{multline*}
and see that this is essentially \eqref{eq:G2metric}, provided one
rescales time $t\mapsto e^{-mt}t$ and  changes names to variables with
the recipe \eqref{eq:coords}.
\end{proof}

For a better understanding of the process one might want to
reconsider it within its symplectic framework
\cite{Hitchin:forms}. The natural variables $p(t),q(t)$ of the candidate
Hamiltonian function $H(t)=H(p,q)$ have to satisfy the standard relations
\begin{displaymath}
      \left\{
        \begin{array} {l}\ds
          p'=-\frac{\partial H}{\partial q}\\[8pt]\ds
           q'=\phantom{-}\frac{\partial H}{\partial p}
        \end{array}
      \right..
\end{displaymath}
They 
translate here into
\begin{displaymath}
    p'=-\tfrac 25 m\sqrt{q+1} \qq q'=\tfrac 65 m\sqrt{p+1},
\end{displaymath}
leading to
%\begin{displaymath}
%    \frac{dp}{dq}=-\frac 26 \sqrt{\frac{q+1}{p+1}},
%\end{displaymath}
%and hence
\begin{equation}\label{eq:levels}
    6\bigl( (p+1)^{3/2}-1\bigr)=1-(q+1)^{3/2}.
\end{equation}
Therefore the functions are
\begin{displaymath}
    p(t)=(1-mt)^{2/5}-1 \text{\q and\q } q(t)=\tfrac 65 m(1-mt)^{1/2}
\end{displaymath}
and $H$, constant on the level curves of \eqref{eq:levels}, is
\begin{displaymath}
    H(t)=\tfrac 25 m(1-mt)^{9/5} + \tfrac 65 m(1-mt)^{3/5}.
\end{displaymath}
\bigbreak

%---------------------------
\subsection{Second example.} 
The Lie algebra
{\small
\bdm
   (-\tfrac 35 m{\rm e}^{17},-\tfrac 35 m{\rm e}^{27},\tfrac
   25 m{\rm e}^{15}-\tfrac 65 m{\rm e}^{37},
   \tfrac 25 m{\rm e}^{25}-\tfrac 65 m{\rm e}^{47},-\tfrac 35 m{\rm e}^{57},
   \tfrac 25 m{\rm e}^{12}-\tfrac 65 m{\rm e}^{67},0),
\edm }
arises from 
$\n=(0,0,\tfrac 25 m\,{\rm e}^{15},\tfrac 25 m\,{\rm e}^{25},0,
\tfrac 25 m\,{\rm e}^{12}) $ and is isomorphic to \eqref{eq:3step}.
Consider the following metric
\begin{multline}
\label{eq:instantonmetric}
    g =
    e^{-\tfrac 45 mt} (dx_1^2 + dx_2^2
+ dx_5^2) +
    \tfrac 9{25}m^2e^{\tfrac 25 m t}(d x_3 + \tfrac 23 x_5
dx_1)^2\\[3pt]
    + \tfrac 9{25} m^2 e^{\tfrac 25 mt}  (d x_4 - \tfrac
23 x_2
    dx_5)^2 + \tfrac 9{25}m^2 e^{\tfrac 25 m t} (d
    x_6 +
\tfrac 23 x_2 dx_1)^2 + e^{-2m t} dt^2
\end{multline}
with the identifications
\begin{displaymath}
  \left\{
  \begin{array} {l}
    {\rm e}^i = e^{\tfrac 35 m t} d x_i,\qquad i = 1, 2, 5,\\[3pt]
    {\rm e}^3 = -\tfrac 15 me^{\tfrac 65 m t}(3dx_3 +2x_5 dx_1),\\[3pt]
    {\rm e}^4 = -\tfrac 15 me^{\tfrac 65 m t}(3dx_4 -2x_2 dx_5),\\[3pt]
    {\rm e}^6 = -\tfrac 15 me^{\tfrac  65 m t}(3dx_6 +2x_2 dx_1),\\[3pt]
    {\rm e}^7 = dt.
  \end{array}
  \right.
\end{displaymath}
These expressions make $g$ a
$G_2$-holonomy metric on $S=\R\times \mathcal{D}$, $\mathcal{D}$ being 
the the $T^3$-bundle over the torus $T^3$ associated to $\text{span}\{{\rm e}_3,
{\rm e}_4, {\rm e}_6 \}$.
\begin{theorem}
The nilpotent Lie algebra $\vs (0,0,\tfrac 25 m{\rm e}^{15},\tfrac 25
m{\rm e}^{25},0,\tfrac 25 m{\rm e}^{12})$ equipped with $\SU(3)$ forms
$\vs \omega_0={\rm e}^{56}-{\rm e}^{23}+{\rm e}^{14},\
\psp_0=-{\rm e}^{345}+{\rm e}^{136}+{\rm e}^{246}+{\rm e}^{125}$ 
generates the Ricci-flat metric
\bdm
  g=(1-mt)^{4/5}g_\text{fibre}+ (1-mt)^{-2/5}g_\text{base}+dt^2
\edm
on $S$, in terms of the flat metrics
$\vs g_\text{fibre}=(\1)^2+(\2)^2+(\5)^2$ and 
$\vs g_\text{base}=(\3)^2+(\4)^2+(\6)^2$ on $\mathcal{D}$.
\end{theorem}
\begin{proof}
The square of $\omega_0$ is an exact form, as 
$d({\rm e}^{364})=\tfrac
12 \oo_0$, whereby
\begin{displaymath}
    \oo(t)=P(t)\,\oo_0
\end{displaymath}
for some smooth function $P$ on an interval $I\set \R$
such that $P(0)=1$. As for $\psp$, the boundary 
conditions ensure that only the term ${\rm e}^{125}$ varies.
Equations \eqref{eq:flow} together
with the primitivity of $\psp$ (holding at all
time) yield
\begin{displaymath}
  \left\{
  \begin{array} {l}
    \psp(t)=-{\rm e}^{345}+{\rm e}^{136}+{\rm e}^{246}+
    (E+1){\rm e}^{125}\\[3pt]
    \omega(t)=\pm\sqrt{P}\omega_0
  \end{array}
  \right.
\end{displaymath}
with $E(0)=1$.
The solution $P(t)=(1-\tfrac 52 t)^{2/5}$ implies that only
the volume of the fibres intervenes in the evolution of the
three-form
\begin{gather*}
  \psp=(1- mt)^{6/5}\,{\rm e}^{125}+{\rm e}^{136}+{\rm e}^{246}-{\rm e}^{345},
\end{gather*}
whilst the K\"ahler structure deforms as
\begin{gather*}
  \omega=(1-mt)^{1/5}\,({\rm e}^{14}-{\rm e}^{23}+{\rm e}^{56}).
\end{gather*}
By writing the compatible metric one obtains the desired expression. 
The fibres grow at double the
speed at which the base shrinks, precisely as in
\eqref{eq:instantonmetric}.\end{proof}

The reader might want to compare the horizontal/vertical split of this metric 
to the similar one of \eqref{eq:BS}.

%%%%%%%%%%%%%%%%%%%%%%%%%%%%
\bibliographystyle{amsplain}

\providecommand{\bysame}{\leavevmode\hbox
  to3em{\hrulefill}\thinspace}
\providecommand{\MR}{\relax\ifhmode\unskip\space\fi MR }
\providecommand{\MRhref}[2]{%
  \href{http://www.ams.org/mathscinet-getitem?mr=#1}{#2}
}
\providecommand{\href}[2]{#2}

\end{document}